\documentclass[final,onefignum,onetabnum]{siamart220329}

\usepackage{braket,amsfonts}

\usepackage{array}

\usepackage[caption=false]{subfig}

\usepackage{pgfplots}

\newsiamthm{claim}{Claim}
\newsiamremark{remark}{Remark}
\newsiamremark{hypothesis}{Hypothesis}
\crefname{hypothesis}{Hypothesis}{Hypotheses}

\usepackage{algorithmic}

\usepackage{graphicx,epstopdf}

\Crefname{ALC@unique}{Line}{Lines}

\usepackage{amsopn}

\usepackage{xspace}
\usepackage{bold-extra}
\usepackage[most]{tcolorbox}

\colorlet{texcscolor}{blue!50!black}
\colorlet{texemcolor}{red!70!black}
\colorlet{texpreamble}{red!70!black}
\colorlet{codebackground}{black!25!white!25}

\lstdefinestyle{siamlatex}{%
  style=tcblatex,
  texcsstyle=*\color{texcscolor},
  texcsstyle=[2]\color{texemcolor},
  keywordstyle=[2]\color{texemcolor},
  moretexcs={cref,Cref,maketitle,mathcal,text,headers,email,url},
}

\tcbset{%
  colframe=black!75!white!75,
  coltitle=white,
  colback=codebackground, 
  colbacklower=white, 
  fonttitle=\bfseries,
  arc=0pt,outer arc=0pt,
  top=1pt,bottom=1pt,left=1mm,right=1mm,middle=1mm,boxsep=1mm,
  leftrule=0.3mm,rightrule=0.3mm,toprule=0.3mm,bottomrule=0.3mm,
  listing options={style=siamlatex}
}

\newtcblisting[use counter=example]{example}[2][]{%
  title={Example~\thetcbcounter: #2},#1}

\newtcbinputlisting[use counter=example]{\examplefile}[3][]{%
  title={Example~\thetcbcounter: #2},listing file={#3},#1}

\DeclareTotalTCBox{\code}{ v O{} }
{ 
  fontupper=\ttfamily\color{black},
  nobeforeafter,
  tcbox raise base,
  colback=codebackground,colframe=white,
  top=0pt,bottom=0pt,left=0mm,right=0mm,
  leftrule=0pt,rightrule=0pt,toprule=0mm,bottomrule=0mm,
  boxsep=0.5mm,
  #2}{#1}

\patchcmd\newpage{\vfil}{}{}{}
\pgfplotsset{compat=1.18}

\usepackage{mathtools}
\usepackage{amsmath}
\DeclareMathAlphabet\mathbfcal{OMS}{cmsy}{b}{n}
\newcommand{\StateDim}{n}

\usepackage[normalem]{ulem}
\usepackage{color}

\newcommand{\js}[1]{{\color{black}{#1}}}

\begin{tcbverbatimwrite}{tmp_output_header.tex}
\title{An Introduction to Solving the Least-Squares Problem in Variational Data Assimilation
\thanks{Submitted to the editors DATE.
\funding{This work was partially funded by a Network Support 
Grant from the Isaac Newton Institute for the Mathematical Sciences and the Engineering \& Physical Sciences Research Council (EPSRC) in the UK (EP/V521929/1).
We also acknowledge funding provided by the Deutsche
Forschungsgemeinschaft (DFG) - Project-ID 318763901 - SFB1294 to visit the University of Potsdam. For the purpose of open access, the authors have applied a Creative Commons Attribution (CC BY) licence to any Author Accepted Manuscript version arising.}}}

\author{I. Dau\v{z}ickait\.{e}\thanks{Parallel Algorithms Team, CERFACS, Toulouse, France (\email{ieva.dauzickaite@cerfacs.fr}, \email{selime.gurol@cerfacs.fr})} \and M.~A. Freitag\thanks{Institute for Mathematics, University of Potsdam, Germany (\email{melina.freitag@uni-potsdam.de}).}
\and S. G\"{u}rol\footnotemark[2]
\and A.~S. Lawless\thanks{School of Mathematical, Physical and Computational Sciences, University of Reading and National Centre for Earth Observation, UK (\email{a.s.lawless@reading.ac.uk})}
\and A. Ramage\thanks{Department of Mathematics and Statistics, University of Strathclyde, UK (\email{a.ramage@strath.ac.uk})}
\and ~~~J.~A. Scott\thanks{School of Mathematical, Physical and Computational Sciences, University of Reading and Scientific Computing Department, Science \& Technology Facilities Council, Rutherford Appleton Laboratory, UK (\email{jennifer.scott@reading.ac.uk})}
\and J.~M. Tabeart\thanks{Department of Mathematics and Computer Science, Eindhoven University of Technology, The Netherlands (\email{j.m.tabeart@tue.nl})}}

\headers{ the Least-Squares Problem in Variational Data Assimilation}
{Dau\v{z}ickait\.{e}, Freitag, G\"{u}rol, Lawless, Ramage, Scott, Tabeart}
\end{tcbverbatimwrite}
\title{An Introduction to Solving the Least-Squares Problem in Variational Data Assimilation
\thanks{Submitted to the editors DATE.
\funding{This work was partially funded by a Network Support
Grant from the Isaac Newton Institute for the Mathematical Sciences and the Engineering \& Physical Sciences Research Council (EPSRC) in the UK (EP/V521929/1).
We also acknowledge funding provided by the Deutsche
Forschungsgemeinschaft (DFG) - Project-ID 318763901 - SFB1294 to visit the University of Potsdam. For the purpose of open access, the authors have applied a Creative Commons Attribution (CC BY) licence to any Author Accepted Manuscript version arising.}}}

\author{I. Dau\v{z}ickait\.{e}\thanks{Parallel Algorithms Team, CERFACS, Toulouse, France (\email{ieva.dauzickaite@cerfacs.fr}, \email{selime.gurol@cerfacs.fr})} \and M.~A. Freitag\thanks{Institute for Mathematics, University of Potsdam, Germany (\email{melina.freitag@uni-potsdam.de}).}
\and S. G\"{u}rol\footnotemark[2]
\and A.~S. Lawless\thanks{School of Mathematical, Physical and Computational Sciences, University of Reading and National Centre for Earth Observation, UK (\email{a.s.lawless@reading.ac.uk})}
\and A. Ramage\thanks{Department of Mathematics and Statistics, University of Strathclyde, UK (\email{a.ramage@strath.ac.uk})}
\and ~~~J.~A. Scott\thanks{School of Mathematical, Physical and Computational Sciences, University of Reading and Scientific Computing Department, Science \& Technology Facilities Council, Rutherford Appleton Laboratory, UK (\email{jennifer.scott@reading.ac.uk})}
\and J.~M. Tabeart\thanks{Department of Mathematics and Computer Science, Eindhoven University of Technology, The Netherlands (\email{j.m.tabeart@tue.nl})}}

\headers{ the Least-Squares Problem in Variational Data Assimilation}
{Dau\v{z}ickait\.{e}, Freitag, G\"{u}rol, Lawless, Ramage, Scott, Tabeart}

\ifpdf
\hypersetup{ pdftitle={An Introduction to Solving the Least-Squares Problem in Variational Data Assimilation} }
\fi

\begin{document}
\maketitle

\begin{tcbverbatimwrite}{tmp_output_abstract.tex}
\begin{abstract}
Variational data assimilation is a technique for combining measured data with dynamical models.  It is a key component of
Earth system state estimation and is commonly used in weather and ocean forecasting. The approach involves a large-scale generalized nonlinear least-squares problem. 
Solving the resulting sequence of sparse linear subproblems requires the use of sophisticated numerical linear algebra methods. 
  In practical applications, the computational demands severely limit the number of iterations of a Krylov subspace solver that can be performed and so
  high-quality preconditioners are vital. 
  In this paper, we present a numerical linear algebra perspective on variational data assimilation and discuss contemporary solution methods for the challenges posed by large-scale geophysical applications. The principal contribution is a focused treatment of the underlying linear algebraic subproblems, accompanied by a concise and clear introduction to the essential concepts of variational data assimilation and an extensive bibliography.

\end{abstract}


\end{tcbverbatimwrite}
\begin{abstract}
Variational data assimilation is a technique for combining measured data with dynamical models.  It is a key component of
Earth system state estimation and is commonly used in weather and ocean forecasting. The approach involves a large-scale generalized nonlinear least-squares problem.
Solving the resulting sequence of sparse linear subproblems requires the use of sophisticated numerical linear algebra methods.
  In practical applications, the computational demands severely limit the number of iterations of a Krylov subspace solver that can be performed and so
  high-quality preconditioners are vital.
  In this paper, we present a numerical linear algebra perspective on variational data assimilation and discuss contemporary solution methods for the challenges posed by large-scale geophysical applications. The principal contribution is a focused treatment of the underlying linear algebraic subproblems, accompanied by a concise and clear introduction to the essential concepts of variational data assimilation and an extensive bibliography.

\end{abstract}




\section{Introduction and motivation} 
Data assimilation is the science of combining information from observations and numerical models to estimate the state of a dynamical system as it evolves over time.  
Although it was initially developed for numerical weather prediction, it is now applied to many classical systems, including geophysical systems such as the Earth’s atmosphere, ocean, and land surface \cite{asch2016data, carrassi_data_2018, Dale91, evensen2022data, freitag_nla_da_2020, FrePot2013, Ghil1991, Kalnay_2002, law2015data, reich2015probabilistic} and, more broadly, to fields such as solar physics \cite{lang2021improving}, ecology \cite{Pandya2022}, cognitive science \cite{Engbert2022}, biology \cite{Evensen2020, lawson1995data} and engineering \cite{moireau2008joint}. 
Here, our focus is on large-scale geophysical systems. 
 
Variational data assimilation (VarDA) corrects the trajectory of the underlying physical dynamical model by incorporating noisy and sparse observations. 
The most probable state of the dynamical system is found by solving an optimization problem. The
cost function of this optimization problem is formulated as a
generalized nonlinear least-squares problem with weights based on
uncertainties in the data and the model.  
VarDA is widely-used in operational weather forecasting. Here, and in the more general context of geophysical systems, VarDA generally displays the following challenging characteristics: large-scale sparse problems, computationally expensive physical models accessible only via operators, and stringent time constraints on obtaining a solution. 

Practical algorithms for tackling the nonlinear least-squares problem are typically based on the truncated Gauss–Newton approach, 
which involves solving a sequence of overdetermined linear least-squares problems of the form
\begin{equation}
\label{Eq:General_form_ls}
\min_{s} \frac{1}{2} \left \| J s - b\, \right \|^2_{W^{-1}}.
\end{equation}
Here and elsewhere, norms $\| s \|_{W^{-1}}^2 = s^T\,W^{-1}s$ are used, where $W$ is symmetric positive definite (SPD). 
The matrix $J$ is rectangular with more rows than columns.
If $J$ is of full rank and the matrix of weights $W$ is SPD then \eqref{Eq:General_form_ls} has a unique solution that can be 
computed approximately using a direct or an iterative method; see the recent book  \cite{Bjorckbook2024} for a comprehensive overview. In VarDA, the computational cost is dominated by the evaluation of the underlying physical
model and so to solve \eqref{Eq:General_form_ls}, a \textit{truncated} iterative method is typically used,
that is, a fixed limited number of iterations of a Krylov subspace solver are computed.
To obtain the required accuracy in the computed solution, it is crucial to accelerate convergence using appropriate \textit{preconditioners} \cite{ brown_multilevel_2016, debreu_multigrid_2016, PearsonPestana2020, tshimanga_limited-memory_2008, wathen_acta}. The design of efficient preconditioners within VarDA raises interesting theoretical and practical questions at the
intersection of optimization and numerical linear algebra.

Our aim in this paper is to provide a unified framework for VarDA from the perspective of numerical linear algebra, with particular emphasis on preconditioning strategies. The paper is aimed primarily at numerical linear algebraists who may be new to VarDA but it will also be of interest to VarDA practitioners seeking a clearer understanding of the associated algebraic structure.
The principal novelty and contribution of this work lies in isolating and clarifying the core linear algebra subproblems that underpin VarDA, presenting them within a coherent and accessible framework. Although the key ideas can be found - explicitly or implicitly - across a wide and scattered body of literature, accessing them typically requires navigating numerous papers written for different audiences. Our treatment brings these elements together succinctly, providing both the necessary VarDA background and a focused discussion of the associated numerical linear algebra challenges. Moreover, we address a significant practical challenge: the language and notation used in the VarDA and numerical linear algebra communities often differ. We seek to provide a coherent and consistent framework that reconciles these conventions, adhering as far as we can to notation and terminology used elsewhere in the literature.
In addition, we include an extensive bibliography.

The remainder of the paper is organised as follows. In Section~\ref{Sec:ProblemStatement}, we introduce the nonlinear least-squares problem, and the associated linear least-squares subproblems (\ref{Eq:General_form_ls}) that arise in VarDA and introduce the terminology and notation that is used throughout the paper. 
The linear least-squares problems can be solved using methods based on 
the normal equations or an augmented system formulation;
this is discussed in Section~\ref{Sec:SolveLinearizedProb}. 
Section~\ref{sec:primaldual} presents the generalized normal equations that arise in VarDA, and considers both the primal and  the dual formulations. Section~\ref{Sec:preconditioning normal equations} looks at 
preconditioning the generalized normal equations.
We discuss  \textit{first-level} preconditioners and then the possibility of accelerating the convergence of the iterative solver further through the employment of a \textit{second-level} preconditioner.
Preconditioning techniques for the augmented system formulation
are presented in Section~\ref{Sec:time-parallel}.
In Section~\ref{sec:Future challenges}, we highlight numerical linear algebra-related challenges that remain in the field of VarDA.
Finally, some concluding remarks are made in Section~\ref{Sec:conclusions}.

\section{The least-squares problems in VarDA}
\label{Sec:ProblemStatement}
Our aim is to estimate the state vector of physical phenomena, such as atmospheric temperature or ocean salinity, in a 
prescribed time interval; this is a common challenge in Earth system modelling. 
Access to the state trajectories over time, $\{x_i \in \mathbb{R}^n \}_{i=0,\ldots,N}$, is obtained through a physical 
dynamical model, $\mathcal{M}_{i}$, which is represented by computationally expensive partial differential equations (PDEs).
Each $\mathcal{M}_{i}$ propagates the state $x_{i-1}$ at time $t_{i-1}$ to the state $x_i$ at time $t_i$ by solving the given 
PDEs. This process includes errors, arising, for example, from a lack of complete representation of physical processes in a numerical model, a lack of model resolution
and numerical discretization errors \cite{howes2017accounting}. These so-called model errors are represented by a time-dependent random variable. Let $q_i$ denote the error in the 
underlying physical model at time $t_i$. We then have
\begin{displaymath}
x_i = \mathcal{M}_{i}(x_{i-1}) + q_{i}, \qquad {i=1,\ldots,N}.
\end{displaymath}
We may also have \textit{a priori} information (known as the background) at the initial time $t_0$, expressed as
\begin{displaymath}
x_b = x_0 + \epsilon_b,
\end{displaymath}
where the background error, $\epsilon_b$, is another random variable. 
We suppose that the state $x_i$ is observed by various instruments, including airborne, ground-based, and space-based sensors. The relation between the observations  $y^o_i \in \mathbb{R}^{m_i}$  and 
the state $x_i$ can be expressed as
\begin{displaymath}
y_i^o = \mathcal{H}_{i}(x_{i}) + \nu_{i}, \qquad {i=0,\ldots,N},
\end{displaymath}
where the observation operator $\mathcal{H}_i$ maps $x_i$ to an $m_i$-dimensional vector representing the state in the 
observation space. This 
operator and the observations themselves again involve errors, which are represented by a time-dependent random variable $\nu_i$, termed the 
observation error. 
The map $\mathcal{H}_i$ may include unit and/or discretization transformations between the state space and the observation space. For instance, the data 
may be observed in radiance, and we are interested in deducing $x_i$, which signifies temperature. Depending on the 
application, $\mathcal{H}_i$ may be complex and nonlinear and, in general, $m_i \ll n$ (the dimension of the state $x_i$).

The goal in data assimilation is to determine the optimal time-distributed state vector 
$x^{\ast} = [(x^\ast_0)^T, (x^\ast_1)^T, \ldots, (x^\ast_N)^T]^T \in \mathbb{R}^{(N+1)n}$ using a given observation set 
$(y_i^o, t_i)$, 
$0 \le i \le N$, the \textit{a priori} state $x_b$ 
and a dynamical model $\mathcal{M}_i$, taking into account their uncertainties. The optimal solution changes with respect to the chosen 
statistical approach and properties of the uncertainties. We assume throughout that the background error 
$\epsilon_b$, the observation error $\nu_i$ and the model error $q_i$ are statistically independent, random variables. In this paper, we consider a Bayesian estimate where 
these errors are assumed to follow a Gaussian distribution with zero mean and SPD covariance matrices 
$B \in \mathbb{R}^{n \times n}$, $R_i \in \mathbb{R}^{m_i \times m_i}$,  and $Q_i \in \mathbb{R}^{n \times n}$, 
respectively. For convenience, we set $p=(N+1)n$.\\

\noindent
\textbf{Weak formulation.} A Bayesian maximum \textit{a posteriori} estimate can be found by solving the following 
\textit{generalized nonlinear least-squares problem} \cite{Stuart_2010}: find 
$x = [(x_0)^T, (x_1)^T, \ldots, (x_N)^T]^T \in\mathbb{R}^p$ that minimizes the quadratic cost function
\begin{equation}
\label{eq:costfunc_weakc_state}  
 \frac{1}{2} \sum_{i=0}^{N} \left\|\mathcal{H}_i(x_i) - y_i^o \right\|_{R_i^{-1}}^2 
+ \frac{1}{2} \left\|x_0 - x_b \right\|_{B^{-1}}^2\\
+  \frac{1}{2}\sum_{i=1}^{N} \left\|x_i - \mathcal{M}_i(x_{i-1}) \right\|_{Q_i^{-1}}^2 .
\end{equation}
This formulation, which
includes the model error, is known in the data assimilation community as  weak-constraint four-dimensional variational data 
assimilation (\textit{weak-constraint 4DVar}) (or the weak formulation) because the model constraint does not need to be satisfied exactly~\cite{sasaki1970some, tremolet_weakconstraint_2006}. \\

\noindent
\textbf{Strong formulation.} If the model error is negligible,  \eqref{eq:costfunc_weakc_state} simplifies to:
find $x_0\in\mathbb{R}^n$ that minimizes
\begin{equation}
\label{eq:costfunc_strong}  
 \frac{1}{2} \sum_{i=0}^{N} \left \|\mathcal{G}_i(x_0) - y_i^o \right \|_{R_i^{-1}}^2 
+ \frac{1}{2} \left \|x_0 - x_b \right \|_{B^{-1}}^2, 
\end{equation}
where $\mathcal{G}_i (x_0) = \mathcal{H}_i (\mathcal{M}_{i}(\cdots \mathcal{M}_{1}(x_0))) = \mathcal{H}_i (x_i) $.
Once the initial state $x_0$ has been determined, the state $x_i$ can be computed using the recurrence relation
$$x_{i} = \mathcal{M}_i(x_{i-1}),\qquad i=1,\ldots,N.$$ 
This form of the model equation acts as a strong constraint on the minimization problem, so it is termed the  \textit{strong constraint 4DVar (or strong formulation)} \cite{sasaki1970some}.

In strong constraint 4DVar, only the initial state $x_0 \in\mathbb{R}^n$ is used in the cost function \eqref{eq:costfunc_strong}. This leads to a simpler formulation compared to weak constraint 4DVar and hence to lower computational costs. Weak constraint 4DVar requires significantly more memory because the cost function \eqref{eq:costfunc_weakc_state} involves the much larger state vector $x \in\mathbb{R}^p$. Furthermore, the model error covariance matrices $Q_i$ are required, although they are not always known. However, the weak constraint formulation is beneficial for longer assimilation windows as it better accounts for the model error build up over time.

In practical applications, the weak and strong formulations\eqref{eq:costfunc_weakc_state} and \eqref{eq:costfunc_strong} have the following important 
properties. 
\begin{itemize}
\item  The dimension $n$ of $x_i$ generally exceeds $10^7$, so the problems are large-scale 
(especially~\eqref{eq:costfunc_weakc_state}) \cite{bauer2015quiet, carrassi_data_2018, chrust2025impact}. 
\item The total number of observations $m = \sum_{i=0}^N m_i$ is small compared to $n$, i.e., $m \ll n$ \cite{bauer2015quiet, 
carrassi_data_2018}.
\item The covariance matrices $B$, $R_i$,  and $Q_i$ ensure the correct weighting is applied to each term in the cost function to obtain the maximum \emph{a posteriori} estimate. They are generally not diagonal and may not be explicitly available. Hence they are modelled, for instance, using diffusion-based correlation operators, or estimated using Monte-Carlo methods, and only their actions on a vector can be computed \cite{bannister2008review, janjic2018representation, tremolet2007model,  weaver2005multivariate}. 
\item $\mathcal{M}_i$  and $\mathcal{H}_i$ are available only as operators. Evaluation of these operators can be computationally expensive (particularly $\mathcal{M}_i$ because it involves solving PDEs). This makes obtaining exact second order derivative information prohibitively expensive \cite{Fisher_PDEconst_2009}.
\end{itemize}
In practice, problems \eqref{eq:costfunc_weakc_state} and \eqref{eq:costfunc_strong} are commonly solved using the truncated Gauss–Newton (TGN) method. TGN is an iterative optimization technique for nonlinear least-squares problems that avoids the explicit computation of the Hessian matrix, which involves second-order derivatives. Instead, the method approximates the Hessian using only first-order derivative (Jacobian) information, resulting in a significant reduction in computational cost for large-scale problems \cite{gratton_approximate_2007, NoceWrig06}.
At each TGN iteration, the nonlinear problem is locally linearized, yielding a linear least-squares problem \eqref{Eq:General_form_ls}. The solution of this problem provides a search direction that is used to update the current estimate. This leads to an inner–outer iteration structure: solving \eqref{Eq:General_form_ls} constitutes the inner iteration, while the repeated linearization and update steps define the outer TGN iteration.

\subsection{Linear least-squares problem for the weak formulation}
When referring to iteration $k$ of the TGN method (that is, the $k$-th outer iteration), we
use the superscript $(k)$. Let the  search direction  from $x^{(k)}$ be $s = [(s_0)^T, (s_1)^T, (s_2)^T, \ldots, (s_N)^T]^T 
\in \mathbb{R}^{p}$. The linear least-squares problem is: find $s$ 
that minimizes the quadratic cost function
\begin{equation} 
\label{Prob:quad_weak}
\begin{split}
      \frac{1}{2} \sum_{i=0}^{N} \|H_i^{(k)}s_i - d_i^{(k)} \|_{R_i^{-1}}^2 
      + \frac{1}{2} &\left \| s_0 - (x_b - x_0^{(k)}) \right \|_{B^{-1}}^2\\
     &+\frac{1}{2} \sum_{i=1}^{N} \left \|s_i - M_i^{(k)} s_{i-1}  - c_i^{(k)} \right \|_{Q_i^{-1}}^2 ,
     \end{split}
\end{equation}
where  $d_i^{(k)} = y_i^o - \mathcal{H}_i(x_i^{(k)}) \in \mathbb{R}^{m_i}$ is known as the innovation,  $c_i^{(k)} = \mathcal{M}_i(x_{i-1}^{(k)}) - 
x_i^{(k)} \in \mathbb{R}^{n}$ is the model error, $H_i^{(k)} \in \mathbb{R}^{m_i \times n}$ is the Jacobian matrix 
of the observation operator $\mathcal{H}_i$ at $x_i^{(k)}$, 
and $M_i^{(k)} \in \mathbb{R}^{n \times n}$ is the Jacobian matrix of the physical model $\mathcal{M}_i$ at $x_{i-1}^{(k)}$. 
Once $s$ is computed, the next iterate is $x^{({k+1})} = x^{({k})} + s$.
This process continues until either the chosen convergence criterion is met or the limit on the allowable 
number of outer iterations is reached. 
In practice, the Jacobian matrices $H_i^{(k)}$ and $M_i^{(k)} \in \mathbb{R}^{n \times n}$ 
 are too large to be stored explicitly. Their action on a vector is commonly computed using automatic differentiation, in which source code that applies the Jacobian to a given vector is generated directly from the source code of the corresponding nonlinear operators \cite{chao1992development, giering1998recipes}.
The term automatic refers to the differentiation process itself and does not necessarily imply the use of automatic compilers. In operational numerical weather prediction systems, such as those at M\'et\'eo France and the European Centre for Medium Range Weather Forecasts, the Jacobian application and its adjoint are implemented through explicitly written, highly optimized manual code (see, e.g., \cite{janiskova1999simplified}).

Dropping the superscript $(k)$ for clarity of notation, \eqref{Prob:quad_weak} can be written in the compact form
\begin{equation}
\label{Prob:compact_weak_quad}  
 \frac{1}{2}\left \|F^{-1} s  - f \right \|_{D^{-1}}^2 
        + \frac{1}{2} \left \|H s - d \right \|_{R^{-1}}^2.
\end{equation}
Here, $d$ is a $m$-dimensional concatenated vector of the $ d_i = y_i^o - \mathcal{G}_i(x_0) \in \mathbb{R}^{m_i}$ and $f = [ (x_b - x_0)^T, (c_1)^T, \ldots, (c_N)^T ]^T\in \mathbb{R}^{p}$. The block rectangular matrix  $H \in \mathbb{R}^{m \times p}$ has the matrices $H_i$ on its main diagonal, i.e., $H = \text{diag}(H_0, H_1, \ldots, H_N)$. Similarly, $R \in \mathbb{R}^{m 
\times m}$ and $D \in \mathbb{R}^{p \times p}$ are SPD block diagonal matrices with $R = \text{diag}(R_0, R_1, \ldots, R_N)$ and $D = \text{diag}
(B, Q_1, \ldots, Q_N)$.
The matrix $F^{-1} \in \mathbb{R}^{p \times p}$ is the block bi-diagonal matrix
\begin{equation}\label{eq:Finv}
F^{-1} =
\begin{bmatrix}
I_n        & 0 & 0 & \cdots & 0 \\
- M_1      & I_n        & 0 & \cdots & 0 \\
0 & - M_2      & I_n        & \ddots & \vdots \\
\vdots     & \vdots     & \ddots     & \ddots & 0 \\
0 & 0 & \cdots     & - M_N  & I_n
\end{bmatrix},
\end{equation}
and its inverse is given by the block lower triangular matrix
\begin{equation}\label{eq:F}
F =
\begin{bmatrix}
I_n        & 0 & 0 & \cdots & 0 \\
M_{1,1}    & I_n        & 0 & \cdots & 0 \\
M_{1,2}    & M_{2,2}    & I_n        & \ddots & \vdots \\
\vdots     & \vdots     & \ddots     & \ddots & 0 \\
M_{1,N}    & M_{2,N}    & \cdots     & M_{N,N} & I_n
\end{bmatrix}.
\end{equation}
where $M_{i,j} = M_j M_{j-1} \ldots M_{i+1} M_i$ represents the sequential application of the Jacobian
matrices of the physical model from $t_{i-1}$ to time $t_j$.
We observe that it is efficient to form matrix-vector products in parallel with $F^{-1}$ (but not with $F$).

The cost function~\eqref{Prob:compact_weak_quad}
is known as the {\textit{weak state formulation}} of 
4DVar~\cite{tremolet_weakconstraint_2006}. 
It can be rewritten as an overdetermined  linear least-squares problem of the form 
(\ref{Eq:General_form_ls}) with
\begin{equation}
\label{eq:overdetermined_leastsquares}
\begin{split}
   J =
   \begin{pmatrix}
       H \\
       F^{-1}
   \end{pmatrix}\in \mathbb{R}^{(m+p) \times p}, \quad b =
      \begin{pmatrix}
       d \\
       f
   \end{pmatrix}\in \mathbb{R}^{m+p}, \\
    W = \begin{pmatrix}
    R & 0 \\
    0 & D
    \end{pmatrix}\in \mathbb{R}^{(m+p) \times (m+p)}.
\end{split}
\end{equation}

\bigskip

\subsection{Linear least-squares problem for the strong formulation} The strong formulation \eqref{eq:costfunc_strong} assumes there is no model error 
and thus it involves only $x_0$. At outer iteration $k$ of TGN, a search direction from $x_0^{(k)}$  is computed by solving the linear least-squares problem 
\begin{equation}
\label{Prob:strong_quad}  
    \min_{s \in \mathbb{R}^n} \left\{
    \frac{1}{2}\left \|G^{(k)} s - d^{(k)} \right \|_{R^{-1}}^2 
        + \frac{1}{2} \left \| s - (x_b- x_0^{(k)}) \right \|_{B^{-1}}^2 \right\}.
\end{equation}
The Jacobian matrix  $G^{(k)} \in \mathbb{R}^{m \times \StateDim}$ represents a concatenation of the $G^{(k)}_i \in \mathbb{R}^{m_i\times \StateDim}$ (the model $\mathcal{G}_i$ linearized about the current iterate $x_0^{(k)}$). As before, $d \in \mathbb{R}^{m}$ and $R = \text{diag}(R_0, R_1, \ldots, R_N)$.
Omitting the superscript $(k)$, this
subproblem can be written as the overdetermined generalized  least-squares problem
(\ref{Eq:General_form_ls}) with
\begin{equation}
\label{eq:overdetermined_leastsquares2}
\begin{split}
   J =
   \begin{pmatrix}
       G \\
       I
   \end{pmatrix} \in \mathbb{R}^{(m+n) \times n}, \quad
   b=    \begin{pmatrix}
       d \\
       x_b - x_0
   \end{pmatrix} \in \mathbb{R}^{m+n},\\
    W = \begin{pmatrix}
    R & 0 \\
    0 & B
    \end{pmatrix} \in \mathbb{R}^{(m+n) \times (m+n)}.
\end{split}
\end{equation}
A summary of the notation that we use for the weak state and strong formulations of the linear least-squares subproblems, 
including the dimensions of the corresponding matrices, is given in Table \ref{tab:notation}. ote that the block matrix $W$ is SPD.

\begin{table}[h]
    \centering
    \begin{tabular}{c|cc|cc}
     &&&&\\
        Notation & Weak state &Dimensions & Strong   & Dimensions\\
        &&&&\\
        \hline
         &&&&\\
        $J$ & $\begin{pmatrix}
            H \\ F^{-1}
        \end{pmatrix}$ & $(m+p)\times p$&$\begin{pmatrix}
            G\\I
        \end{pmatrix}$ & $(m+n)\times n$\\[0.5cm]
        $b$ & $\begin{pmatrix}
            d\\f
        \end{pmatrix}$ & $m+p$& $\begin{pmatrix}
            d\\x_b-x_0
        \end{pmatrix}$ & $m+n$\\[0.5cm]
        $W$  &$\begin{pmatrix}
            R & 0 \\ 0 & D
        \end{pmatrix}$ & $(m+p)\times (m+p)$& $\begin{pmatrix}
            R & 0 \\ 0 & B
        \end{pmatrix}$ & $(m+n)\times(m+n)$ \\
    \end{tabular}
    \caption{Summary of the notation for the weak state and strong formulations of the 
    linear least-squares subproblems.}
    \label{tab:notation}
\end{table}

\section{Solving linear least-squares problems}
\label{Sec:SolveLinearizedProb}
Having introduced the linear least-squares problems
corresponding to the weak state and strong formulations, in
this section, we discuss methods for solving general large-scale overdetermined  least-squares problems of the form (\ref{Eq:General_form_ls}).
There are a number of possible methods, see, e.g., the recent
book \cite{Bjorckbook2024} and the survey article \cite{sctu_acta:2025} for comprehensive discussions.
Here we consider two commonly-used approaches that are particularly relevant to VarDA.
\\

\noindent
\textbf{Normal equations.} Solving the linear least-squares problem \eqref{Eq:General_form_ls} is mathematically equivalent to solving the 
\emph{generalized normal equations}  given by
\begin{equation}
\label{eq:normal_equations}
( J^T W^{-1} J )\, s =  J^T W^{-1}\, b .
\end{equation}
The SPD {weighted normal matrix} $J^T W^{-1} J$ is the Hessian of the quadratic problem~\eqref{Eq:General_form_ls}. \\

\noindent
\textbf{Augmented system.} Problem \eqref{Eq:General_form_ls} can be reformulated as a constrained optimization 
problem~\cite{fisher2018low,fisher_parallelization_2017}, for which the Karush-Kuhn-Tucker (KKT) conditions represent a 
special case of the \emph{augmented system}~\cite{benb:1999,Bjorckbook2024, orban2017iterative}
\begin{equation}
\label{eq:augmented system}
K \begin{pmatrix}
\lambda \\
s 
\end{pmatrix}
= 
\begin{pmatrix}
b \\
0 
\end{pmatrix}, \quad \mbox{with} \quad K = 
\begin{pmatrix}
 W &J \\
 J^T & 0
\end{pmatrix}.
\end{equation}
Here $\lambda\in \mathbb{R}^{m+p}$ is a vector of Lagrange multipliers. The matrix $K$ is a sparse symmetric indefinite matrix 
that
is non-singular if $J$ is of full rank and $\mathcal{N}(W) \, \cap \, \mathcal{N}(J^T) = \{ 0\}$.  
Saddle-point problems of this form arise in a wide variety of practical problems~\cite{benzi_numerical_2005}. 
In contrast to many applications, 
in VarDA the $(2,1)$ block (the so-called constraint block) is much more expensive to apply than the $(1,1)$ block,
although  matrix-vector products with $J$ (involving $F^{-1}$ in the weak state formulation) can be implemented in parallel 
\cite{fisher2018low,fisher_parallelization_2017}. \\

\noindent
There are advantages and disadvantages associated with both the normal equations and the augmented system. Solving an SPD system is often regarded as more straightforward than solving an indefinite system and there is a wider class of applicable solution methods (see below for a brief introduction). However, the augmented system
can have favourable numerical properties \cite{Arioli1987}. 
For example, the condition number quantifies the sensitivity of a problem to  perturbations to the data. A matrix with a large condition number
is said to be {ill-conditioned}, otherwise it is well-conditioned. For the rectangular matrix $W^{-1/2}J$,  the 2-norm condition number $\kappa(W^{-1/2}J)$ is equal to the ratio of its largest to its smallest nonzero singular values. Forming the
normal equations squares the condition
number, that is, $\kappa(J^TW^{-1}J) = \kappa^2(W^{-1/2}J)$.
\js{We note that the conditioning of the augmented system matrix $K$ can be controlled by considering a scaled matrix $P_L K P_R$, where $P_L = \text{diag}(W^{-1/2},\alpha^{-1} I)$ and $P_R = \text{diag}(\alpha W^{-1/2}, I)$. The optimal value of the parameter $\alpha$ ensures $\kappa(P_LKP_R)\leq 2\kappa(W^{-1/2}J)$. Computing the optimal $\alpha$
is expensive as it involves the smallest nonzero singular value of $W^{-1/2}J$ but heuristic approaches are available \cite{Arioli1987}.}

Other potential problems with using the normal equations include the 
loss of information when the inner products used to compute 
the entries of the normal matrix are accumulated (see the discussion in \cite{Bjorckbook2024}).
Even if the inner products are accumulated in double precision
arithmetic, a serious loss of
information can occur when the computed normal matrix is stored in the working precision. In general, whenever $\kappa(W^{-1/2}J) \ge \epsilon^{-1/2}$ (where $\epsilon$ is the machine precision) we can
expect the computed normal matrix to be singular (or indefinite).
Furthermore, solution methods that explicitly form 
and factorize the normal equations are not backward
stable because it can be shown that the best backward error bound
contains a factor $\kappa(W^{-1/2}J)$ \cite{high:02}.
Nevertheless, using the normal equations is attractive, particularly if only modest
accuracy is required or if the problem can be preconditioned to make it well conditioned; for VarDA, we discuss preconditioning in Section~\ref{Sec:preconditioning normal equations}.

In VarDA, the augmented system \eqref{eq:augmented system} is significantly larger than the 
normal equations and solution methods can be prohibitively expensive in terms of the memory requirements. 
However, $K$ is sparse whereas the normal matrix can be much  denser (for instance, if $J$ contains
a single dense row then the normal matrix is dense). Note that, by eliminating $\lambda$ in~\eqref{eq:augmented system}, we 
recover~\eqref{eq:normal_equations}. 

Equations \eqref{eq:normal_equations}  and
\eqref{eq:augmented system} are examples of large-scale linear systems of equations in which the system matrix is symmetric and the right-hand side vector is known. 
There are many methods for solving such systems; see, for example, the books 
\cite{Bjorckbook2024,liesen_strakos_2012, orban2017iterative, saad2003,  sctu_book:2023},  the review article \cite{dars:16}, and the references therein. 
The methods can be split into two main classes:
direct  and iterative (with hybrid methods combining techniques from both classes). We briefly discuss these in the following subsections.

\subsection{Direct methods}
Direct methods use a finite sequence of elementary transformations
to rewrite the system matrix as a product of simpler matrices
in such a way that solving linear systems with these matrices is 
relatively straightforward. Provided $J$ is of full rank and $W$ is SPD,
a Cholesky factorization of the symmetrically permuted normal matrix 
\[\Pi^TJ^T \,W^{-1}J\,\Pi=LL^T\]
can be computed, where the (square) 
factor $L$ is a lower triangular matrix and the permutation matrix $\Pi$ is chosen to preserve sparsity in $L$.
For the symmetric indefinite augmented system \eqref{eq:augmented system}, we can compute a  factorization
\[\Pi_K^T \,K\,\Pi_K = L_K \,\Delta_K \, L_K^T,\]
where $L_K$ is a 
unit lower triangular matrix and $\Delta_K$ is block diagonal
with diagonal blocks of size 1 and 2. 
In this case, the permutation matrix $\Pi_K$ is chosen
to retain sparsity and for numerical  stability.
The need to ensure stability makes factorizing sparse symmetric indefinite matrices
much harder than sparse SPD matrices, necessitating the use of  sophisticated algorithms
and implementations. These become even more sophisticated if parallel implementations are sought.

Having computed the factorization,  linear
systems with the triangular factors can be 
solved using simple forward and back substitutions. Unfortunately, it is challenging to obtain good
speedups for this solve step  in a parallel environment because the substitution steps are inherently serial, although there has been work
on circumventing this, for example  by using Jacobi iterations \cite{chow2018}.

An alternative approach for solving linear least-squares problems 
that avoids forming the normal matrix
or the augmented system is to compute a QR factorization (see, e.g., \cite{Golub2012}). This
 seeks to express  $W^{-1/2}J$ as a product of an orthogonal matrix
and an upper triangular matrix.
While this can offer greater
numerical stability, it is a more expensive approach (in terms of time and memory).

 Implementing sparse direct algorithms so that the resulting software is
 efficient and robust is complicated,
 requiring significant experience and expertise. However, when applied appropriately, direct methods can  provide black-box solvers for computing  solutions with predictable accuracy.
  The main shortcomings of direct methods are that they can require a large number of 
  numerical operations and a large amount of memory. These
  demands increase with the size of the system matrix and its density,
  and eventually become prohibitive. A further limitation is that direct methods require explicit access to the system matrix; as a result, they are unsuitable for linear systems where access is only indirectly available  through vector products with an operator, as is the case for VarDA.
  Moreover, they can compute solutions to an accuracy that may not be either needed
  or warranted by the supplied data.   Consequently, for the very large systems that arise in
  data assimilation,  iterative methods are used.

\subsection{Preconditioned iterative methods}

Iterative methods for solving a generic (square) linear system of equations of the form
\[\mathbb{A}w=b\] 
aim to compute a sequence of  
approximate solutions that converges to the required solution in an acceptable 
number of iterations. The most commonly-used methods are Krylov subspace methods \cite{saad2003, vdV2003}. 
The system matrix $\mathbb{A}$ does not need to be stored explicitly as it is only used indirectly, through matrix--vector products.
How much storage is required depends on the iterative method and on whether it is necessary to incorporate 
reorthogonalization between some (or all) of the vectors generated. For  methods where the orthogonal vectors can be 
calculated using a short-term recurrence relation, such as the well-known conjugate gradient method (CG) for SPD systems 
\cite{hest:52} and MINRES for general symmetric linear systems \cite{pasa:75}, in theory only a small number of vectors of 
length the size of the linear system need
to be stored.
However, in finite precision arithmetic, there can be a loss of orthogonality that can adversely affect the rate of convergence. It may therefore be advantageous  to keep (some of) the previously computed
vectors and employ  reorthogonalization \cite[\S7.5]{Demmel1997}. Other popular iterative methods, including 
GMRES \cite{saad1986}, have no short-term recurrence and the number of vectors 
that must be held and the computational costs increase with the 
iteration count. In this case, it may be necessary to include strategies (such as restarting) 
to limit the work and storage needed.

An advantage of iterative solvers is that the user can choose how many iterations to perform or specify the required accuracy in the computed solution.
Properties that influence the rate of convergence are the initial solution guess, the right-hand side vector, and the 
system matrix $\mathbb{A}$.
The conditioning of $\mathbb{A}$ is of particular importance.
For a square matrix $\mathbb{A}$ of full rank the
2-norm condition number is $\kappa(\mathbb{A}) = \|\mathbb{A}\|_2 \|\mathbb{A}^{-1}\|_2$. If $\mathbb{A}$ is SPD, this becomes $\kappa(\mathbb{A})=\lambda_{max}(\mathbb{A})/\lambda_{min}(\mathbb{A})$, where $\lambda_{max}$ and
$\lambda_{min}$ are the largest and smallest eigenvalues of $\mathbb{A}$.

To illustrate the importance of the conditioning on the performance of an iterative solver,
it can be shown that if $\mathbb{A}$ is SPD then the 
approximate solution $w^j$ at iteration $j$ of  
the CG method satisfies the bound 
\begin{displaymath} 
\|w-w^j\|_{\mathbb{A}}\le 2\left(\frac{\sqrt{\kappa(\mathbb{A})}-1}{\sqrt{\kappa(\mathbb{A})}+1}\right)^j 
\|w-w^0\|_{\mathbb{A}} .
\end{displaymath} 
However, this error bound can be highly pessimistic. In particular, it does not show the potential for the CG method to converge
superlinearly, or that the rate of convergence depends on the distribution of all the eigenvalues of 
$\mathbb{A}$. For example, if the eigenvalues of $\mathbb{A}$ lie in clusters, the convergence of CG can be significantly faster than the above bound implies, \js{see, for example \cite[\S5.1]{NoceWrig06}, \cite[\S2]{kelley1995}, \cite[\S3]{greenbaum1997} and \cite{liesen_strakos_2012}}.

Because the normal matrix $J^T W^{-1}J$ in \eqref{eq:normal_equations} is SPD,  an obvious approach is to employ the 
CG method or its Lanczos variant (Lanczos-CG)~\cite{Frommeretal99,pasa:75,saad2003} to solve the normal equations. 
The CGLS method for linear least-squares problems is derived by a slight algebraic rearrangement of the CG method  
\cite{hest:52}. 
This involves additional storage and work per iteration but has the advantage that the least-squares residual 
is recurred, rather than the residual of the normal equations; this is discussed in \cite{Bjorck1998}.

The well-known LSQR method \cite{pasa:82} is another Lanczos-type algorithm for solving least-squares problems.
It is again mathematically equivalent to CG applied
to the normal equations but can offer improved numerical stability, especially when the system matrix is ill-conditioned and many iterations are needed to achieve the requested accuracy.
Applying the Lanczos process to the augmented system~\eqref{eq:augmented system} with $W=I$ forms the basis of the 
Golub-Kahan lower bidiagonalization procedure used in LSQR~\cite{benb:1999}. Thus, the least-squares 
problem~\eqref{Eq:General_form_ls} can be solved either by employing a change of variables
or by using the generalized G-LSQR approach~\cite{benb:1999,orban2017iterative}, which
is based on a generalized Golub-Kahan bidiagonalization technique.
A potential disadvantage of the LSQR and G-LSQR methods is that they require additional storage compared to CG; for least-squares problems \eqref{eq:overdetermined_leastsquares} and \eqref{eq:overdetermined_leastsquares2}, the extra storage is equal to the total number of observations $m$.

The LSMR method \cite{fosa:2011} is also based on Golub-Kahan bidiagonalization. It is mathematically equivalent to the MINRES method applied to the normal equations,
with  both the least-squares residual and the normal equations residual decreasing  monotonically. This may allow LSMR to terminate after fewer iterations
than CGLS and LSQR \cite{gosc:2017}.

In practical VarDA applications, CG is commonly used, without explicitly forming the  potentially ill-conditioned normal matrix.
Incorporating reorthogonalization has been found to be crucial for solution 
accuracy \cite{el_akkraoui_preconditioning_2013, Guroletal_2014}.
The large-scale nature of the problems 
and the prohibitive cost of applying the system matrix mean that only a small number of iterations are performed. This truncation of the solver makes reorthogonalization feasible because 
it is  only requires the storage of a
corresponding number of vectors. 

 The CG method is generally not used for solving the augmented system \eqref{eq:augmented system} because, for indefinite systems,
there is no guarantee that it will not fail. Hence other Krylov
subspace methods are employed, including MINRES or GMRES or quasi-minimal residual (QMR) methods \cite{freund91},
such as SQMR \cite{freund94}. 
When CG is applied to the normal equations \eqref{eq:normal_equations},  
the energy norm of the error decreases monotonically, thereby implicitly minimizing the cost function of the least-squares problem \eqref{Eq:General_form_ls} over the Krylov subspace built during the CG iterations~\cite{Golub2012}. 
For \eqref{eq:augmented system}, however, an important consideration  is that, although for MINRES and GMRES the augmented 
system residual decreases monotonically, because $s$  forms only part of the 
solution vector,  this cost function  can increase as the iteration
count increases \cite{dauzickaite2022preconditioning, gratton_guaranteeing_2018}. This is of particular concern in VarDA where truncated Krylov subspace methods are standard and thus, when  the iterations are stopped, the value of the  cost function can be larger than at the start. 
The safe-guarded method proposed in \cite{gratton_guaranteeing_2018} ensures sufficiently many iterations are performed for each inner solve to obtain an improved solution to \eqref{eq:costfunc_weakc_state}. This comes at the expense of additional evaluations of the cost function. 

Preconditioning aims to speed up convergence of an iterative method by transforming the given system
into an equivalent system (or one from which it is easy to recover the solution
of the original system) that has `nicer' numerical properties. 
Conceptually, this involves replacing the original system by the modified equations
\begin{displaymath}
 P^{-1}\mathbb{A}w=P^{-1}b, \qquad\mbox{or}\qquad
 \mathbb{A}P^{-1}\hat{w}=b,\quad w=P^{-1}\hat{w},
\end{displaymath}
where $P$ is the \textit{preconditioner}.
These represent so-called {left} and {right}
preconditioning. 
If $P$ is SPD, it can also be applied symmetrically via a
factorization; this is termed split preconditioning. In all three cases, it is
necessary only to solve systems with  $P$,
without explicitly computing $P^{-1}$ or its factors.
$P$ should be chosen such that 
the conditioning of the preconditioned problem is better than that of the original problem, ideally with a more favorable 
eigenvalue distribution and it should be inexpensive (i.e., the cost of its construction and application should be
less than the resulting savings in the iterative solver runtime). Preconditioners can be easily incorporated into Krylov subspace methods leading to, e.g., the 
well-known preconditioned conjugate gradient (PCG) algorithm \cite{Concus_1976}. In practice, however, determining a good 
preconditioner is highly problem dependent and can be very challenging.
In VarDA, the early  truncation of the solver after a fixed number of iterations means that the role of preconditioning in accelerating convergence in the initial iterations is particularly important.
Preconditioners that have been used in this context will be discussed in Sections 5 and 6.

\section{A note on primal and dual formulations}
\label{sec:primaldual}
For the strong and weak state linear least-squares problems introduced in Section~\ref{Sec:ProblemStatement}, 
using (\ref{eq:overdetermined_leastsquares2})
and (\ref{eq:overdetermined_leastsquares}) respectively, the normal 
equations (\ref{eq:normal_equations}) take the following (primal) forms:
\begin{equation}
\label{eq:strong_primal}
 \underbrace{\left(B^{-1} + G^T R^{-1} G \right)}_{\mathbb{A}_S} s = \underbrace{B^{-1}(x_b - x_0) + G^TR^{-1}d}_{b_S} ,
\end{equation}
and
 \begin{equation}
\label{eq:weak_state_primal}
   \underbrace{\left(F^{-T}D^{-1}F^{-1} + H^T R^{-1} H \right)}_{\mathbb{A}_W} s = \underbrace{F^{-T}D^{-1} f + H^TR^{-1}d}_{b_W}.
\end{equation}
Introducing the change of variables $v = F^{-1} s$,  \eqref{eq:weak_state_primal} can also be considered using the so-called 
\textit{weak forcing formulation}~\cite{fisher_parallelization_2017}
\begin{equation}
\label{eq:weak_forcing_primal}
  \underbrace{\left(D^{-1} + F^T H^T R^{-1} H F \right)}_{\mathbb{A}_F} v = \underbrace{D^{-1}f + F^TH^TR^{-1}d}_{b_F}. 
\end{equation}
Note that $\mathbb{A}_S\in\mathbb{R}^{n\times n}$ and $\mathbb{A}_W, \mathbb{A}_F\in\mathbb{R}^{p\times p}$.
When $m \ll n$ (i.e., when there are far fewer observations than
the dimension of the state space), the computational cost and memory requirements for solving 
these primal problems (including the memory needed for reorthogonalization) can be reduced by applying iterative methods to the associated \textit{dual problems} in 
$m$-dimensional space. 
For the strong formulation \eqref{eq:strong_primal}, the dual problem involves solving the $m \times m$
linear system
\begin{equation}
\label{eq:dual_linear_system}
   \left(R^{-1}G B G^T + I \right) u = R^{-1} \left( d - G (x_b - x_0)\right),
\end{equation}
and then computing $s =  x_b - x_0 + B G^T u$ \cite{Courtier97, GrattonTshimanga2009}. 
The physical space statistical analysis system (PSAS)~\cite{CohndaSiGuo98} was the 
first dual approach to be proposed, using the $R$-inner product~\cite{gratton_gurol_toint_2013}. Conventional implementations 
of PSAS with diagonal $R$ employ $R^{-1}$ as a preconditioner via the square-root 
$R^{-1/2}$~\cite{Courtier97,AkkraouiGauthier2010,Akkraouietal2008},
and solve
\begin{equation}
\label{eq:PSAS}
    \left(R^{-1/2}G B G^TR^{-1/2} + I \right) z = R^{-1/2} \left( d - G (x_b - x_0)\right), \quad u = R^{-1/2}z,
\end{equation}
using CG or Lanczos-CG furnished with the canonical inner product. However, as illustrated in~\cite{AkkraouiGauthier2010, 
GrattonTshimanga2009},  the dual iterates of PSAS produce a non-monotonic decrease of the cost function 
in the linear least-squares problem~\eqref{Prob:strong_quad}. 
One possible remedy is
to use  MINRES to solve \eqref{eq:PSAS}. Numerical results in \cite{AkkraouiGauthier2010} illustrate that, in this case, the cost function decreases monotonically. 

Another possibility is to use the restricted preconditioned
conjugate gradient (RPCG) method~\cite{GrattonTshimanga2009}, which also
solves~\eqref{eq:dual_linear_system} using a CG method but equipped with the
(possibly semi-deﬁnite) $GBG^T$-inner product instead of the $R$-inner product. RPCG generates, in exact arithmetic, 
the same iterates as PCG applied to~\eqref{eq:strong_primal} with preconditioner $B$
(under assumptions on the initial guess and preconditioner choice that are easily satisfied).
The Lanczos variant of RPCG is introduced in~\cite{Guroletal_2014}. 

The matrix
$R^{-1}G B G^T + I$ in~\eqref{eq:dual_linear_system} is nonsymmetric but the system can be solved using a non-standard 
inner product within CG or Lanczos-CG~\cite{bramble1988,gratton_gurol_toint_2013,gurol_phdthesis_2013, StollWathen_2008}.
If the iterative method is run to full convergence 
then the computed solution obtained from the dual problem is mathematically equivalent to the solution of the linear 
system~\eqref{eq:strong_primal}. However, this equivalence is not guaranteed if CG is truncated early and neither is the 
monotonicity of the cost function evaluated using the dual iterates. 

The systems (\ref{eq:weak_state_primal}) and (\ref{eq:weak_forcing_primal}) arising from the
weak formulations can be solved using a dual approach in analogous ways. In this case, the computational gains are even 
more significant than for the strong formulation, as the dimension of the problem to be solved again reduces to the number of 
observations.

\section{Preconditioning the normal equations}
\label{Sec:preconditioning normal equations}

In VarDA, it is common to theoretically transform the Hessian of the optimization problem to a new operator with a more 
favorable eigenvalue distribution  using so-called \textit{first-level preconditioning}; this is discussed in 
Section~\ref{section:Firstlevelprecond}.
In Section~\ref{Sec:SecondLevelPrecond} 
we explore improving convergence further by  preconditioning the preconditioned normal equations; this is termed \textit{second-level preconditioning}.
Note that  combining preconditioners is widely used in other
fields; see, for instance, \cite{al2021, Tangetal2009}. Finally, in Section~\ref{sec:non-constant A}, we consider (second-level) preconditioners for a sequence of slowly varying linear systems, as encountered in practice in full VarDA simulations.
We denote by $A$ the matrix ${\mathbb A}$ {\bf after the application of the first-level preconditioner.}

\subsection{First-level preconditioning}
\label{section:Firstlevelprecond}

A good choice for a first-level preconditioner depends on the problem characteristics, which may relate to its 
physical properties or the algebraic structure of the resulting linear system \cite{Benzi2002, benzi_numerical_2005, 
MalekStrakos2014, PearsonPestana2020, wathen_acta}. 

Consider the normal matrix $\mathbb{A}_S$ in (\ref{eq:strong_primal}).
Because the number of observations is much smaller than the size of the 
state vector, the term $G^TR^{-1}G$ is a low-rank update of 
 $B^{-1}$. The matrix $B$  is often highly ill-conditioned, leading to ill-conditioned normal equations \cite{haben_conditioning_2011-1, 
haben_conditioning_2011, lorenc1997development, tabeart2018conditioning}.
In VarDA applications, the most common first-level preconditioning step applies a {split preconditioner}
through a change of variables. Using a factorization 
$B=UU^T$
leads to the symmetric preconditioned system
\[ U^T \mathbb{A}_S\, U z  =  U^T b_S , \quad s = Uz.\] 
If $U$ is square and of full rank, it acts as a perfect scaling,  that is, the components of the transformed variable $z$ are mutually 
uncorrelated with unit variance \cite{menetrier_overlooked_2015}. 
It is important to note that  $U$  need not be obtained via a Cholesky factorization; in practice,  a suitable $z$ is chosen from physical principles and then $U$ is defined as a series of linear transformations from $z$ to $s$, enabling its application without explicit matrix construction~\cite{Derber1999,fisher2003,weaver_courtier_2001}.
The preconditioned normal matrix becomes 
\[
{A}_S = I + U^TG^TR^{-1}G\,U,
\] 
where the second term  has rank $m\ll n$ and is a low-rank
update of $I$.
${A}_S$ has $n-m$ eigenvalues clustered at one and the remainder are greater than one. This transformation of the spectrum is
expected to improve the convergence of Krylov subspace methods
\cite{haben2011_MO, NoceWrig06, tabeart_new_2022}. 

When the factor $U$ cannot easily be estimated, the PCG method without the explicit use of a factorization can be 
employed \cite{destouches2024qg,el_akkraoui_preconditioning_2013}. Matrix-vector products with $B^{-1}$ can be avoided by 
introducing an auxiliary vector~\cite{derber_Bprecond_89}. Alternatively, 
$\mathbb{A}_S$ can be preconditioned by $B$ from the left or right, leading to the matrices $I + 
BG^TR^{-1}G$ and $I + G^TR^{-1}GB$, respectively. These have the same eigenvalue spectrum as ${A_S}$ 
\cite{el_akkraoui_preconditioning_2013}. However, because  symmetry is not preserved, standard CG methods cannot be 
applied. A bi-conjugate gradient (Bi-CG) method is used in \cite{el_akkraoui_preconditioning_2013},
which shows that the PCG algorithm introduced in \cite{derber_Bprecond_89} is a particular case of Bi-CG applied to the data assimilation problem. 
In \cite{gurol_phdthesis_2013, Guroletal_2014}, it is noted that  the CG algorithm can be adapted to solve the preconditioned problem through the use of a non-standard inner product~\cite{bramble1988, Liesen_Strakos_2008, StollWathen_2008, saad2003}. In exact arithmetic,  this 
produces iterates that are mathematically equivalent to those 
obtained with a split preconditioner.

For the weak formulations, the normal matrices are  $\mathbb{A}_W$ in \eqref{eq:weak_state_primal} and $\mathbb{A}_F$ in \eqref{eq:weak_forcing_primal}. 
These again comprise a full-rank term  plus a low-rank update. Their condition numbers have different sensitivities to the parameters of the assimilation process 
and neither is consistently superior \cite{el-said_conditioning_2015}. The matrix $D = \text{diag}
(B, Q_1, \ldots, Q_N)$ in the full-rank term includes the error covariance matrices $B$ and $Q_i$ and is often highly ill-conditioned. It is therefore natural for any preconditioning strategy to treat this term first. For the weak forcing formulation \eqref{eq:weak_forcing_primal}, the  structure of  $\mathbb A_F$ is similar to that
of $\mathbb A_S$ \eqref{eq:strong_primal}. In practice, the $Q_i$ are constructed such that a factorization is available and hence, a first-level split-preconditioner can be based on a factorization of the form $D=D_1 D_1^T$. The preconditioned normal matrix is then
\[
    {A}_F = I + D_1^TF^TH^TR^{-1}HFD_1. \]
Observe that matrix-vector products with ${A}_F$ require computationally expensive products with $F$ that cannot readily be parallelized.

For the weak state  formulation,  a preconditioner of the form $\tilde{F}^{-T}D^{-1}\tilde{F}^{-1}$ can be used. In VarDA, $\tilde{F}^{-1}$ is typically constructed by replacing $M_i$ in \eqref{eq:Finv} by an approximation $\tilde{M}_i$ such that products with $\tilde{F}$ can be performed in parallel. For example,  $\tilde{M}_i = 0 $ and  $\tilde{M}_i = I$ \cite{gratton_note_2018}.
In this case, the condition number of the preconditioned matrix $(\tilde{F}D\tilde{F}^T)(F^{-T}D^{-1}F^{-1})$ is bounded.
An alternative is to select $\tilde{M}_i =  \tilde{M}$,
where $\tilde{M}$ seeks to incorporate  information from the model. The resulting preconditioner can be applied by exploiting the resulting Kronecker structure of $\tilde{F}$ \cite{palitta_stein-based_2023}. 
Note that even if a matrix $\tilde{F}^{-1}$ is a good approximation to $F^{-1}$, the matrix $\tilde{F}^{-T}D^{-1}\tilde{F}^{-1}$  can be an 
arbitrarily poor approximation to $F^{-T}D^{-1}F^{-1}$
\cite{brpe:86,gratton_note_2018,wathen_comments_2022}.

Another approach is to note that $FD^{1/2}=D^{1/2}+L$, where $L$ is a strictly lower triangular matrix.
A  preconditioner of the form $D^{1/2}+\tilde{L}$ can be defined by taking $\tilde{L}$ to be a low-rank approximation to 
$L$. One possibility is to use randomized methods \cite{dauzickaite_timeparallel_2021}. Although these require expensive 
matrix-vector products with $M_{i,j}$, $\tilde{L}$ is returned as a truncated singular value decomposition and is thus cheap 
to apply.

\subsection{Second-level preconditioning}
\label{Sec:SecondLevelPrecond}
First-level preconditioners
cluster many of the eigenvalues at one but some large eigenvalues may remain and 
these can hinder the convergence of the iterative solver. 
To accelerate convergence, a second-level preconditioner may be needed. 

One way of preconditioning $A$ is to approximate its inverse. 
In VarDA, this is typically done by using limited memory preconditioners (LMPs)  as
second-level preconditioners \cite{gratton2011class,tshimanga_limited-memory_2008}.
LMPs are defined by
\begin{equation}
\label{eq:lmp}
P = [I-Z(Z^TAZ)^{-1}Z^TA] [I-AZ(Z^TAZ)^{-1}Z^T] + \theta Z(Z^TAZ)^{-1} Z^T,
\end{equation}
where  $Z$ is a matrix with $\ell$ linearly independent columns, and $\theta > 0$ is a scaling parameter that is often set to $1$~\cite{ Fisher_PDEconst_2009, GiraudGratton2006}.
Note that if $Z$ spans the entire finite-dimensional space and $\theta = 1$ then $P = A^{-1}$.
This family of LMPs is inspired by the BFGS method \cite{NoceWrig06} for minimizing a nonlinear cost function by gradually approximating the inverse of the Hessian.

Special cases of the LMP arise for particular choices of the columns of $Z$. Let the eigenpairs of $A$ be $(z_i, \lambda_i)$ 
with the $z_i$ orthonormal and $\lambda_1 \ge \lambda_2 \ge \ldots \ge \lambda_\ell > 1$.  If $Z=
[z_1,\ldots,z_\ell]$ then the so-called {spectral LMP}~\cite{gratton2011class} or deflating preconditioner~\cite{FrankVuik2001, GiraudGratton2006} is given by
\begin{equation}
\label{eq:spectralLMP}
 P_{\text{spec}} = I+ Z (\theta\Lambda^{-1}-I) Z^T = I-\sum_{i=1}^\ell \left(1-\frac{\theta}{\lambda_i}\right)z_iz_i^T,
\end{equation}
where $\Lambda=\text{diag}(\lambda_1,\ldots,\lambda_\ell)$. Note that the factorization 
$P_{\text{spec}} = P_{\text{spec}}^{1/2}P_{\text{spec}}^{1/2}$ can easily be obtained by replacing $\lambda_i$ and $\theta$ in \eqref{eq:spectralLMP} with their square roots \cite{Fisher1998,tshimanga_limited-memory_2008}.

When applied to $A$ with $\theta = 1$, the LMP preconditioner
(\ref{eq:lmp}) adds at least $\ell$ eigenvalues to the cluster at one, while the rest of the spectrum does not
expand~\cite{gratton2011class}. 
Other values of $\theta$ lead to different positions of the eigenvalue cluster. In \cite{Mouhtaletal2026},  
selecting different values is proposed so that $ P_{\text{spec}}$ is not only a good approximation for $A^{-1}$ but also 
effectively reduces the energy norm of the error (which CG monotonically minimizes), particularly in the early iterations. 
Connections with the deflated CG method are also established, offering further insights into selecting $\theta$ at a 
negligible cost. 
The importance of selecting an appropriate scaling parameter is also highlighted in \cite{gratton2011class}.
Choosing $\theta$ to minimize the condition number of the preconditioned matrix is suggested.

 LMPs  have similarities with other preconditioning techniques in the literature. In \cite{tshimanga2007} it is shown that applying CG preconditioned by \eqref{eq:lmp}, with $\theta=1$ and a specific choice of initial point 
 is analytically equivalent to the deflated CG method \cite{Saadetal2000} with the columns of $Z$ forming the deflation subspace. 
In domain decomposition, this LMP corresponds to the balancing Neumann-Neumann approach (BNN)~\cite{gratton2011class, Mandel1993}. The connection between a two-level multigrid operator, BNN, and deflation methods is established in \cite{NabbenVuik2006, Tangetal2010, Tangetal2009}.

\subsection{Second-level preconditioning for a sequence of systems}
\label{sec:non-constant A}

In a full VarDA simulation, a sequence of slowly varying
linear systems must be solved.  For the second-level
preconditioners described in Section~\ref{Sec:SecondLevelPrecond}, 
the normal matrix $A$ at each outer iteration $k$ was treated
as being fixed (that is, we assumed in constructing the preconditioner that the system matrices did not change with $k$). In this section,
we consider the case that these matrices are not fixed, but vary (usually by a small amount) at each outer iteration. So
after the application of first-level preconditioning, at outer iteration $k$ we need to solve a linear system of equations with system matrix $A^{(k)}$. Typically,
we want to exploit information generated when solving system $k$ 
in preconditioning system $k+1$, or use some other prior knowledge from a previous iteration. Success in constructing good preconditioners depends on the matrices $A^{(k)}$  not 
changing rapidly with $k$.

Approximations of the dominant eigenvalues and corresponding eigenvectors of $A^{(k)}$  can be obtained using the Lanczos method or computed within the PCG iteration itself \cite{Golub2012}. These approximate eigenpairs are known as Ritz pairs; they can be used within  \eqref{eq:lmp} to precondition system $k+1$. This is  the Ritz LMP, 
whilst using Ritz pairs within \eqref{eq:spectralLMP} is termed the inexact spectral LMP~\cite{gratton_approximate_2007}. The latter is employed in  operational weather forecasting, where only converged Ritz pairs are used \cite{Fisher1998, Fisher_PDEconst_2009, Saadetal2000}. Perturbation analysis is presented in \cite{GiraudGratton2006}. 
When applied with converged Ritz pairs, the inexact spectral LMP exhibits similar behavior to the Ritz LMP, although the latter necessitates storing one additional vector. 

Quasi-Newton LMPs \cite{gratton2011class, MorNoc2000} 
choose the columns of $Z$ from the search directions of PCG. 
When all available search
directions or Ritz vectors are used, the quasi-Newton LMP and the Ritz LMP are mathematically equivalent in exact arithmetic. However, the former has twice the storage cost.

When using information coming from PCG or \js{Lanczos-CG}, an LMP can only be used on the second and subsequent outer iterations. Second level preconditioning of the initial system ($k=1$) remains an issue.
This is considered  in \cite{gratton_reduced_2011}.
In the numerical linear algebra literature, there has been
significant emphasis on using randomized algorithms to approximate the eigenspectrum of SPD matrices~\cite{Frangellaetal2023, Halkoetal2011}.  These ideas have been employed to approximate the inverse matrix on each outer iteration (including the first)~\cite{dauzickaite_timeparallel_2021,dauzickaite_randomised_2021,scotto_phdthesis_2022,subrahmanya_randomized_2024}. 
This approach has the additional advantage of being applicable even if $A^{(k)}$ varies significantly with $k$.

Some elements of multigrid and multilevel solvers have been used for preconditioning VarDA problems.
In \cite{debreu_multigrid_2016}, a multigrid V-cycle is applied as a preconditioner 
for $A^{(k)}$ at each outer iteration. A multilevel limited memory approximation to the inverse of $A^{(k)}$ (based on eigenvalue
decompositions obtained from several coarser grid levels) has also been used as a second-level preconditioner \cite{brown_multilevel_2016}.
 When used in conjunction with a local Hessian decomposition, this can result in savings of both
computational time and memory compared to the standard spectral LMP preconditioner (\ref{eq:spectralLMP}).

For the dual formulation, the LMP can be constructed to be symmetric with respect to the $GBG^T$-inner product~\cite{gurol_phdthesis_2013}, consistent with the inner product used in RPCG~\cite{gratton_gurol_toint_2013}.
Second-level preconditioning is also used in~\cite{Egbert1997,souopgui_comparison_2017}, where the importance of employing 
different inner products is emphasized. Symmetry with respect to the inner product 
must be maintained throughout the outer iterations when using PCG or its Lanczos equivalent. In the dual formulation, when the LMP is constructed using the Ritz pairs obtained from the previous outer loop, symmetry is not necessarily preserved. A strategy  proposed in~\cite{gratton_gurol_toint_2013} ensures global convergence through a trust-region approach. Alternatively,~\cite{scotto_phdthesis_2022} uses randomized algorithms based on an 
inner product that inherently preserves symmetry.

\section{Preconditioning the augmented system}
\label{Sec:time-parallel}
In this section, we discuss preconditioning approaches for the generalized augmented system formulation
\eqref{eq:augmented system}. Although the strong formulation \eqref{eq:strong_primal} can  be written as an augmented system 
\cite{rao_time-parallel_2016}, to date most work on preconditioning 
has focused on the weak formulation. 
In this case,  \eqref{eq:augmented system} can be expressed as a  $3\times 3$ block saddle-point system
\begin{equation}\label{eq:augmentedWC}
K \begin{pmatrix}
\lambda \\
s 
\end{pmatrix} =
    \begin{pmatrix}
        R & 0 & H\\
        0 & D& F^{-1} \\
        H^T & F^{-T} & 0
    \end{pmatrix}\begin{pmatrix}
        \lambda_o \\ \lambda_b \\ s
    \end{pmatrix}= \begin{pmatrix}
        d\\f\\0
    \end{pmatrix}.
\end{equation}
 Both the normal equations \eqref{eq:weak_state_primal} and the augmented system \eqref{eq:augmentedWC} can be described as `time-parallel' or `all-at-once' \cite{gander201550} because 
 matrix-vector products with the system matrix 
 only involve  $F^{-1}$ and $F^{-T}$
 (and not $F$ or $F^T$), avoiding sequential products with the $M_{i,j}$ operators (recall \eqref{eq:F}). 
 
 On modern computer architectures, real-time speed-ups can be achieved by distributing operations with $M_i$ and $M_{i,j}$ 
 over many processors. The study \cite{moore2023weak} reports that, although solving \eqref{eq:augmentedWC} requires more 
inner iterations than solving \eqref{eq:weak_state_primal} to achieve a comparable reduction in the cost function, exploiting 
time-parallel algorithms can  reduce the total computational time. Further  improvements can potentially be achieved by 
solving the linear system in a lower precision  than the precision used in the outer iteration, and by using a lower 
resolution linearized model \cite{moore2023weak}.
 
 There is a wealth of research in the numerical linear algebra literature devoted to preconditioners for Krylov subspace methods for saddle-point systems; see, for instance, the survey articles \cite{benzi_numerical_2005,PearsonPestana2020,rees_optimal_2010}
and the references therein. The most successful approaches exploit the block structure of
$K$, possibly together with physical information about the blocks. In VarDA, the aim is to design preconditioners that are time-parallel, taking into account  the cost  of applying the different blocks within $K$ (and their inverses).

The (negative) Schur complement of $K$ with respect to
$W=\begin{pmatrix}
            R & 0 \\ 0 & D
        \end{pmatrix}$
is $S = F^{-T}D^{-1}F^{-1} + H^TR^{-1}H$ (which is the normal
matrix $\mathbb A_W$ in \eqref{eq:weak_state_primal}). The basic block diagonal  preconditioner and its inverse 
are  given by 
\[P_D = \begin{pmatrix}
    R &0&0\\ 0& D&0 \\ 0& 0& S 
\end{pmatrix},  \quad P_D^{-1} = \begin{pmatrix}
    {R}^{-1}&0&0 \\0& D^{-1}&0 \\ 0&0& {S}^{-1} 
\end{pmatrix}.\]

\smallskip  
Block triangular preconditioners \cite{bramble1988} are of the form

\smallskip
\[P_T = \begin{pmatrix}
    R &0 &  H\\0 & D & F^{-1}  \\0 &0 & {S} 
\end{pmatrix}, \quad P_T^{-1} = \begin{pmatrix}
    {R}^{-1} & 0 & -R^{-1}HS^{-1} \\ 0 & D^{-1} & -{D}^{-1}F^{-1}{S}^{-1}\\ 0& 0& {S}^{-1} 
\end{pmatrix}. \] 

The cost of applying $P_T^{-1}$ is higher than for $P_D^{-1}$
because it involves an additional multiplication by $J=\begin{pmatrix}
            H \\ F^{-1}
        \end{pmatrix}$.
For both $P_D$ and $P_T$ , $S^{-1}$ is typically replaced by a computationally affordable approximation $\tilde{S}^{-1}$. 
This can be obtained using the methods developed for preconditioning the normal equations (see Section 
\ref{Sec:preconditioning normal equations}). %

Preconditioners that approximate $J$, referred to as inexact constraint preconditioners, have been well studied  \cite{bergamaschi_inexact_2007,bergamaschi_erratum_2011,sesanaSpectralAnalysisInexact2013}. This motivated the development of a data assimilation-specific preconditioner in which 
         $F$ is approximated by $\tilde{F}$, giving
         
\[P_C = \begin{pmatrix}
    R &0 &0 \\ 0& D& \tilde{F}^{-1} \\  0 & \tilde{F}^{-T}  &0
\end{pmatrix},  \quad 
P_C^{-1} = \begin{pmatrix}
    {R}^{-1}& 0&0\\0& 0 & \tilde{F}^T \\0&\tilde{F} & -\tilde{F}D\tilde{F}^T\end{pmatrix}.\]
Key advantages of $P_C$ are that setting the (1,3) block to zero
greatly simplifies the computation of $P_C^{-1}$ and the application of $D^{-1}$ (which may  
not be available as an operator) is avoided. This preconditioner has been reported to reduce the  iteration count compared to 
the block diagonal and block triangular preconditioners \cite{freitag_low-rank_2018, green_model_2019, tabeart_saddle_2023}. 
In a similar spirit to the second-level preconditioners for the normal equations, 
 information from previous outer iterations can be used to update the preconditioners.
This has been applied to $P_C$ to find a low-rank update to the approximation $\begin{pmatrix}
    0 \\ \tilde{F}^{-1}
\end{pmatrix}$
to $J$  and its transpose \cite{fisher2018low}. 

Much  work  has  focused on developing computationally feasible approximations $\tilde{F}$ of $F$ that can be used within 
$P_D$, $P_T$ and $P_C$. Many of these replace $M_i$ in \eqref{eq:Finv} with some approximation $\tilde{M}_i$ to facilitate 
parallel computation 
\cite{fisher_parallelization_2017,freitag_low-rank_2018,Gratton2016_weakconstraint,gratton_guaranteeing_2018, 
tabeart_saddle_2023,palitta_stein-based_2023}. 

Observe that the $3 \times 3 $ augmented system \eqref{eq:augmentedWC} can be reduced to a $2\times 2$ saddle-point problem 
in which the $(2,2)$ block is nonzero \cite{dauzickaite2022preconditioning,dauzickaite_spectral_2020}.  However, preliminary 
numerical explorations indicate that this  system can  suffer from non-monotonicity of the linear least-squares problem cost function, and slow 
convergence. Currently, in VarDA there is a lack of preconditioners  for this reduced form.

\section{Future challenges}\label{sec:Future challenges}

From the discussions above, it can be seen that in VarDA many  interesting challenges relating to numerical linear algebra remain.
Here, we briefly summarize some of these.\\

\begin{itemize}
\item Operationally, the linear systems in VarDA must be solved using an iterative solver. However, because of the 
computational costs (in terms of time and possibly also memory) the solver is not run ``to convergence'', but is terminated 
after a fixed (typically small) number of iterations. Hence, we need to understand how classical asymptotic results for 
iterative solvers apply in the context of early stopping. 
For the augmented system formulation, the non-monotonicity of the cost function of the linear least-squares problem in the early iterations is particularly problematic
and needs to be addressed when combined with preconditioning.
 \item Preconditioning the linear systems is a huge challenge. Current first-level preconditioning strategies are very standard and, as far as we are aware, the only second-level preconditioners
employed in practice are LMPs. For the augmented system approach,  the constraint block is expensive to apply and so many  standard preconditioning techniques are not applicable. There is scope for exploring more sophisticated preconditioning techniques, in particular methods that are tailored towards the specific (physical) application and model problem used in VarDA. 
 Furthermore, more advanced preconditioners are needed that 
 seek to exploit recent developments in the data assimilation system  \cite{destouches2024mlblue, destouches2024qg,goux_impact_2023,goux_impact_2025,tabeart_new_2022}.
 More work also needs to be targeted at the dual formulation, which is potentially attractive when $m \ll n$.

\item Randomization has been considered within preconditioning strategies for VarDA and for replacing the iterative solver entirely \cite{bousserez_enhanced_2020}, but practical algorithms for large-scale problems have yet to be developed. Randomized algorithms may also be beneficial for speeding up other computational tasks when
solving the structured least-squares problems described in this paper.
\item Machine learning has only been used fairly recently \cite{ackmann2020machinelearning} in preconditioning for linear systems. The cyclic nature of VarDA and the availability of data may allow machine learning strategies to inform the design of preconditioners for VarDA \cite{trappler_state-dependent_2025}.

\item Data assimilation problems and the corresponding least-squares problems are becoming ever larger. Adapting current methods is challenging and will require the 
exploitation of new hardware and modern parallel architectures. Mixed precision algorithms have the potential to deliver improved performance and might be particularly suitable in the limited budget setting. Experimental results on employing lower precision for the linearized model in VarDA show that stabilization techniques are essential for Krylov subspace methods even when medium-complexity models are used \cite{hatfield2020single}; moving to large-scale models and applying reduced precision in other components of the process brings additional challenges.

\item When the \textit{a priori} error, observation error and model error are assumed to be independent Gaussian random variables, the Bayesian estimate results in nonlinear weighted least-squares problems (as in \eqref{eq:costfunc_weakc_state} and \eqref{eq:costfunc_strong}) where the weights, given by the SPD covariance matrices, define energy norms, i.e., weighted Euclidean norms. When we relax this Gaussian assumption for some or all of the errors, the Bayesian inference problem no longer results in a concise minimization problem of the form 
\eqref{eq:costfunc_weakc_state} or \eqref{eq:costfunc_strong}. In the most general cases sampling methods like MCMC will have to be used to find the posterior distribution \cite{Bardsley2018}. 

\end{itemize}

Finally, we observe that, in this paper, we have focused on the challenge of solving least-squares problems in VarDA. We have not covered other questions and issues that also require sophisticated numerical linear algebra techniques within sequential and variational data assimilation, for example, Kalman filtering, low-rank approximations, or dimension and model reduction. Further details are given in \cite{freitag_nla_da_2020}. 

\section{Concluding remarks}
\label{Sec:conclusions}
The main goal of this paper is to introduce the key concepts of variational data assimilation to the numerical linear algebra community, using a unified framework with consistent terminology and notation to 
summarise a wide range of concepts and ideas. In particular, we have shown how variational data assimilation requires the solution of a sequence of large sparse linear least-squares problems with a specific structure.  
We have summarized the main points in the solution of those least-squares problems using preconditioned iterative methods applied to the normal equations and the augmented system formulation. 
The focus in VarDA, in particular for large-scale geophysical systems, is on the first few iterations of the iterative solver because the computational costs involved in each iteration 
and the possible time restrictions in practical applications mean that early stopping is typical. This is a major distinction compared to many other areas in which iterative solvers are run to convergence and means that preconditioning of the linear systems is essential
to improve accuracy. We have presented an overview of the approaches to preconditioning that are employed when solving the linear systems arising in VarDA. In addition, we have provided a literature review that will be particularly useful to those in
the numerical linear algebra community
who are unfamiliar with the field of variational data assimilation but who are motivated  to investigate and engage with some of the many remaining challenges that we  highlight in Section~\ref{sec:Future challenges}. \js{For more comprehensive treatments of data assimilation, we refer the interested reader to the books~\cite{asch2016data,Kalnay_2002,reich2015probabilistic}.}

It is important to appreciate that for large scale problems implementing the methods discussed in this paper is extremely challenging.
As a result, there is currently limited  availability of software for testing new ideas. There have been some steps towards this with generic software packages being made available, including DAPPER: Data Assimilation with Python: a Package for Experimental Research \cite{Raanes2024}. This is aimed only at testing ideas on small models.
Recent efforts within the operational weather forecasting community have focused on developing more accessible frameworks to simplify the implementation of numerical algorithms and strengthen collaboration with research institutions. Examples include OOPS (Object-Oriented Prediction System) and JEDI (Joint Effort for Data Assimilation Integration) \cite{tremolet2020joint}, both of which have publicly available versions.

\bigskip
\paragraph{Acknowledgments}
We would like to thank our colleague Nancy Nichols for commenting on a draft of this manuscript and joining with us in 
useful discussions. We are also grateful to the anonymous reviewers for their helpful and constructive comments.

\bibliographystyle{siamplain}
\bibliography{DA}

\end{document}